\documentclass[reqno,12pt]{amsart}

\RequirePackage{amssymb}
\RequirePackage{times}
\RequirePackage{graphicx}
\RequirePackage{hyperref}
\RequirePackage{lipsum}
\RequirePackage[normalem]{ulem}
\RequirePackage[hang,small,bf]{caption}

\newcommand{\R}{\mathbb{R}}
\newcommand{\C}{\mathbb{C}}

\newcommand{\E}{\mathbb{E}}

\newcommand{\bP}{\mathbb{P}}

\def\Ab{\mathrm{Ab}}
\def\Tr{\mathrm{Tr}}

\newtheorem{theorem}{Theorem}
\newtheorem{lemma}[theorem]{Lemma}

\newtheorem{corollary}[theorem]{Corollary}

\makeatletter
\g@addto@macro{\endabstract}{\@setabstract}
\newcommand{\authorfootnotes}{\renewcommand\thefootnote{\@fnsymbol\c@footnote}}%
\makeatother

\makeatletter
\let\@fnsymbol\@arabic
\makeatother

\begin{document}
\begin{center}
  \LARGE 
  On the Markus-Spielman-Srivastava inequality for sums of rank-one matrices 
	\par \bigskip
  \normalsize
  \authorfootnotes
  Vladislav Kargin\footnote{Email: vladislav.kargin@gmail.com; current address: 282 Mosher Way, Palo Alto, CA 94304, USA} \par \bigskip

\end{center}

\begin{center}
\textbf{Abstract}
\end{center}

\begin{quotation}
We extend the result of Markus, Spielman, and Srivastava about the sum of rank-one symmetric random matrices to the case when the isotropy assumption on the random matrices is relaxed.   
\end{quotation}

Consider the sum $\sum_{i=1}^{m} v_i v_i^{\ast}$, where $v_i \in \C^d$ are independent random vectors and $v_i$ are bounded in expected norm. What can be said about the norm of this sum? Several concentration bounds were given by Rudelson \cite{rudelson99}, Ahlswede and Winter \cite{ahlswede_winter02}, and others (\cite{tropp15}). These bounds show that with high probability the norm of the sum is less than a threshold logarithmic in $d$. Recently,  Markus, Spielman, and Srivastava have shown that the norm of the sum is less than a threshold independent of $d$, with a positive probability. Their result is formulated for isotropic vectors, that is, when $\E\sum_{i=1}^{m} v_i v_i^{\ast}=I$. With some modifications, their methods are applicable to the non-isotropic case as well. In particular, in this paper we prove the following result.

\begin{theorem}
\label{main}
Let $v_1$, \ldots, $v_m$ be independent random vectors in $\C^d$ with finite support. Suppose that  $\E\sum_{i=1}^{m} v_i v_i^{\ast} \leq I_d$ and $\E\|v_i\|^2\leq \varepsilon$ for all $i$. Then 
\begin{equation}
\bP\left[\|\sum_{i=1}^{m} v_i v_i^{\ast}\|\leq (1+\sqrt{\varepsilon})^2 \right]>0.
\end{equation}
\end{theorem}

For the case  $\E\sum_{i=1}^{m} v_i v_i^{\ast} = I_d$, this inequality was proved in \cite{mss13}. Our theorem shows that the assumption of isotropy is irrelevant. What is essential is that $\|\E\sum_{i=1}^{m} v_i v_i^{\ast}\| \leq 1$.

 Markus, Spielman, and Srivastava used a corollary of their inequality to prove the Kadison-Singer conjecture. Below is the analogue of the corollary for the non-isotropic case. 

\begin{corollary}
\label{corollary_partition}
Let $r$ be a positive integer and let $u_i\in \C^d$ be $m$ vectors such that $\sum_{i=1}^m u_i u_i^{\ast} \leq I$ and
 $\|u_i\| \leq \delta$. Then, there exists a partition $\{S_1,\ldots,S_r\}$ of $\{1,\ldots,m\}$ such that for all $k$, 
\begin{equation}
\|\sum_{i\in S_k} u_i u_i^{\ast}\| \leq \frac{1}{r}(1+\sqrt{r\delta})^2
\end{equation} 
\end{corollary}

\textbf{Proof of Corollary \ref{corollary_partition}:} Define vectors $w_{i,j}\in \C^{rd}$: 
\begin{equation}
w_{i,1}=
\left ( 
\begin{array}{l}
u_i \\
0^d \\
\vdots \\
0^d \\
\end{array}
\right),
w_{i,2}=
\left ( 
\begin{array}{l}
0^d \\
u_i \\
\vdots \\
0^d \\
\end{array}
\right),
\ldots,
w_{i,r}=
\left ( 
\begin{array}{l}
0^d \\
0^d \\
\vdots \\
u^i \\
\end{array}
\right).
\end{equation}
Let $v_i$ be random vectors that take value $\sqrt{r} w_{i,k}$ with probability $1/r$.
Then $\E \|v_i\|^2 = r\|u_i\|^2 \leq r\delta$ and $\sum_{i=1}^{m}\E v_i v_i^{\ast}\leq I_{rd}$. Then by Theorem \ref{main}, there exists an assignment $v_i=\sqrt{r} w_{i,k}$, such that 
\begin{equation}
\left\|\sum_{k=1}^{r} \sum_{i:v_i=w_{i,k}} (\sqrt{r} w_{i,k})(\sqrt{r} w_{i,k})^{\ast}\right \| \leq (1+\sqrt{r\delta})^2.
\end{equation}
Set $S_k=\{i:v_i=w_{i,k}\}$, then for all $k$, 
\begin{equation}
 \left\| \sum_{i\in S_k} u_i u_i^{\ast}\right \| \leq \frac{1}{r} (1+\sqrt{r\delta})^2.   
\end{equation}
\hfill $\square$

\textbf{Proof of Theorem \ref{main}:} We follow the lines of proof in \cite{mss13}. The key ingredient is the following formula: 
\begin{equation}
\label{formula_det}
\E\det\left(xI-\sum_{i=1}^{m} v_i v_i^{\ast}\right)
=\prod_{i=1}^{m}(1-\partial_{z_i})\left.\det\left(xI+\sum_{i=1}^{m} z_i\E v_i v_i^{\ast})\right)\right|_{z_1=\ldots=z_m=0},
\end{equation}
which is proved as Theorem 4.1 in \cite{mss13}.
Recall that a multivariate polynomial $p(z_1,\ldots,z_n)$ is called \emph{real stable} if all its coefficients are real and if $p(z_1,\ldots,z_n) \ne 0$, whenever $\Im(z_i)>0$ for all $i$.
The determinant on the right-hand side of (\ref{formula_det}) is a real stable polynomial in variables $x$ and $z_i$ by fundamental Proposition 2.4 in \cite{borcea_branden08}. It follows from the properties of real stable polynomials that the polynomial on the left-hand side of \ref{formula_det} is also real stable. Since it is a polynomial in one variable, hence all its roots are real. Moreover, the convex combinations of real stable polynomials have important root interlacing properties. These properties imply (as in Theorem 4.5 in \cite{mss13}) that there are such values $w_i$ of random variables $v_i$ that the largest root of  $\det\left(xI-\sum_{i=1}^{m} w_i w_i^{\ast}\right)$ is at least as large as the largest root of $\E\det\left(xI-\sum_{i=1}^{m} v_i v_i^{\ast}\right)$. Hence, the conclusion of Theorem \ref{main} follows from Theorem \ref{main2} below. \hfill $\square$

\begin{theorem}
\label{main2}
Let $A_1$, \ldots, $A_m$ be Hermitian positive  semidefinite matrices. Suppose that  $\sum_{i=1}^{m} A_i \leq I$ and $\Tr A_i \leq \varepsilon$ for all $i$. Then the largest root of the polynomial 
\begin{equation}
\prod_{i=1}^{m}(1-\partial_{z_i})\left.\det\left(xI+\sum_{i=1}^{m} z_i A_i)\right)\right|_{z_1=\ldots=z_m=0}
\end{equation}
is at least as large as $(1+\sqrt{\varepsilon})^2$.
\end{theorem}

In the proof  we will use the following notation. We say that $x\in \R^m$ is above the roots of a multivariate polynomial $p(z_1,\ldots,z_m)$ if $p(z_1,\ldots,z_m)>0$
 for all $z\geq x$ (that is, if $z_i\geq x_i$ for all $i$). The set of all points above the roots of $p$ is denoted $\Ab_p$. 

Let $p(z_1,\ldots,z_m)$ be a real stable polynomial and $z \in \Ab_p$. Then \textit{the barrier function} of $p$ at point $z$ 
in direction $i$ is defined as 
\begin{equation}
\Phi_p^i(z)=\partial_{z_i}\log p(z)=\sum_{j=1}^r \frac{1}{z_i-\lambda_j(z_1,\ldots,z_{i-1},z_{i+1},\ldots,z_m)},
\end{equation}
where $\lambda_j(z_1,\ldots,z_{i-1},z_{i+1},\ldots,z_m)$ are roots of the polynomial
 $q_i(t)=p(z_1,\ldots,z_{i-1},t,z_{i+1},\ldots,z_m)$.

We will also use the following lemma (Lemmas 5.9 and 5.10 in \cite{mss13}).
\begin{lemma}
\label{lemma_barrier_shift}
Suppose that $p(z_1,\ldots,z_m)$ is real stable, that $z \in \Ab_p$, and that 
\begin{equation} 
\delta\geq \frac{1}{1-\Phi_p^i(z)}>0.
\end{equation}
Then  for all $i$ and $j$,  $z \in \Ab_{(1-\partial_{z_j})p}$ and 
\begin{equation} 
\Phi_{(1-\partial_{z_j})p}^i(z+\delta e^j)\leq \Phi_p^i(z).
\end{equation}
\end{lemma}

\textbf{Proof of Theorem \ref{main2}:} Let $P(x,z_1,\ldots,z_m)=\det(xI+\sum_{i=1}^{m}z_i A_i)$.
Let $x=(1+\sqrt{\varepsilon})^2$ and $t=-1-\sqrt{\varepsilon}$. Then, 
$x+t>0$. It follows (since $\sum_{i=1}^{m} A_i\leq I$) that $xI+t\sum_{i=1}^{m} A_i$ is positive definite. This implies that $(x,z)\in \Ab_P$ for all $z \geq t1$, where $1$ denote the vector $(1,\ldots,1)$. 

Next, 
\begin{eqnarray}
\Phi_P^i(x,z)&=&\partial_{z_i} \log \det(xI+\sum_{i=1}^{m}z_i A_i)\\
             &=&\Tr \partial_{z_i} \log (xI+\sum_{i=1}^{m}z_i A_i) \\
						 &=& \Tr \left[(xI+\sum_{i=1}^{m}z_i A_i)^{-1}A_i\right], \label{inequality_for_barrier}
\end{eqnarray}
where the last equality follows by differentiating the power expansion for the logarithmic function and using the fundamental property of trace ($\Tr(XY)=\Tr(YX)$).
Since $\sum_{i=1}^{m} A_i \leq I$ and all $A_i$ are positive semidefinite, hence for all $z\geq t1$ we have $xI+\sum_{i=1}^{m}z_i A_i\geq (x+t)I$ and $(xI+\sum_{i=1}^{m}z_i A_i)^{-1}\leq (x+t)^{-1}I$. 

\begin{lemma}
\label{trace_product}
Let $X$ and $Y$ be Hermitian and positive semidefinite. Then $\Tr(XY)\leq \|X\|\Tr(Y)$.
\end{lemma}
Proof of Lemma \ref{trace_product}: Note that $\Tr(XY)=\Tr(Y^{1/2}X^{1/2}X^{1/2}Y^{1/2})=\| X^{1/2}Y^{1/2}\|^2_2$, where $\|A\|_2$ denotes the Frobenius norm of $A$. Then the claim of the lemma follows from the inequality $\|AB\|_2 \leq \|A\| \|B\|_2$. (Section 5.6, Exercise 20 on page 313 in \cite{horn_johnson85}).

\hfill $\square$

 By Lemma \ref{trace_product}, we infer from (\ref{inequality_for_barrier}) that 
\begin{equation}
\Phi_P^i(x,z)\leq (x+t)^{-1}\varepsilon,
\end{equation}
for all $z \geq t1$.

Hence, if $\delta:=-t$, then 
\begin{equation}
\frac{1}{1-\Phi_P^i(x,z)}\leq \frac{1}{1-(x+t)^{-1}\varepsilon}=\delta
\end{equation}
for all $z \geq t1$.
By starting with $(x,t1)$ and applying Lemma \ref{lemma_barrier_shift} sequentially to each of the components of vector $t1$, we find that $x$ is above the roots of polynomial 
\begin{equation}
\label{mixed_char_polynomial}
\prod_{i=1}^{m}(1-\partial_{z_i})\left.\det\left(xI+\sum_{i=1}^{m} z_i A_i)\right)\right|_{z_1=\ldots=z_m=0},
\end{equation}
which is equivalent to the statement of the theorem. 
 \hfill $\square$

\bibliographystyle{plain}
\bibliography{comtest}

\end{document}